\documentclass[12pt]{article}
\usepackage{amsfonts}
\usepackage{latexsym}
\usepackage{amsmath}
\usepackage{amssymb}
\usepackage{amssymb}
\usepackage{amsthm}
\usepackage{mathrsfs}

\newtheorem{thrm}{Theorem}[section]

\newtheorem{lemma}[thrm]{Lemma}

\newtheorem{prop}[thrm]{Proposition}

\newtheorem{cor}[thrm]{Corollary}

\newcommand{\newsection}{    
\setcounter{equation}{0}\section}
\def\appendix#1{\addtocounter{section}{1}\setcounter{equation}{0}
\renewcommand{\thesection}{\Alph{section}}
\section*{Appendix \thesection\protect\indent \parbox[t]{11.15cm}{#1}}
\addcontentsline{toc}{section}{Appendix \thesection\ \ \ #1}}

\newcommand{\be}{\begin{eqnarray}}
\newcommand{\ee}{\end{eqnarray}}
\newcommand{\bea}{\begin{eqnarray}}
\newcommand{\eea}{\end{eqnarray}}
\newcommand{\ba}{\begin{array}}
\newcommand{\ea}{\end{array}}

\newenvironment{frcseries}{\fontfamily{frc}\selectfont}{}

\def \la {\label}

\def\a{\alpha}

\font\mybb=msbm10 at 11pt

\def\bb#1{\hbox{\mybb#1}}

\def\bZ {\bb{Z}}

\def\bC {\bb{C}}

\def\CP{\bC{\mathrm P}}

\begin{document}
\begin{titlepage}
\begin{center}
\vspace*{-1.0cm}


\vspace{2.0cm} {\Large \bf Vanishing theorems  on $(\ell|k)$-strong
 K\"ahler manifolds
with torsion
} \\[.2cm]

\vspace{1.5cm}
 {\large S.~ Ivanov${}^1$ and  G. Papadopoulos$^2$}

\vspace{0.5cm}

${}^1$ University of Sofia, Faculty \\of Mathematics and
Informatics,\\ blvd. James Bourchier 5, \\1164, Sofia, Bulgaria

\vspace{0.5cm}
${}^2$ Department of Mathematics\\
King's College London\\
Strand\\
London WC2R 2LS, UK\\

\vspace{0.5cm}

\end{center}

\vskip 1.5 cm
\begin{abstract}
We derive sufficient  conditions for the vanishing of plurigenera,
$p_m(J), m>0$,  on compact $(\ell\vert k)$-strong,
$\omega^\ell\wedge
\partial\bar\partial \omega^k=0$,  K\"ahler manifolds with
torsion. In particular, we show that the plurigenera of closed
$(\ell\vert k)$-strong manifolds, $k<n-1$, for which ${\rm
hol}(\hat\nabla)\subseteq SU(n)$ vanish, where $\hat\nabla$ is
the Hermitian connection with skew-symmetric torsion. As a
consequence all generalized k-Gauduchon manifolds for which ${\rm
hol}(\hat\nabla)\subseteq SU(n)$ do not admit holomorphic (n,0)
forms. Furthermore we show that all conformally balanced,
$(\ell\vert k)$-strong K\"ahler manifolds with torsion,
$k\not=n-1$,  are  K\"ahler.
 We also give several examples of $(\ell\vert k)$-strong K\"ahler and Calabi-Yau manifolds with torsion.

\end{abstract}

\end{titlepage}


\setcounter{section}{0}
\setcounter{subsection}{0}



\newsection{Introduction}

Hermitian manifolds have widespread applications in both physics
and differential geometry. These are complex manifolds equipped
with a metric $g$,  $g(JX, JY)=g(X,Y)$, and a Hermitian form
$\omega(X,Y)=g(X, JY)$ which  is  (1,1) with respect to the
complex structure $J$. There are many examples of Hermitian
manifolds as every complex manifold admits a Hermitian structure.
In many applications, Hermitian manifolds have  additional
properties which are expressed as either a condition  on $\omega$
or as a restriction on the holonomy of one of the Hermitian
connections. A condition on  the Hermitian form of a
2n-dimensional  manifold is
\bea
\omega^\ell\wedge
\partial\bar\partial \omega^k=0~, ~~~1\leq k+\ell\leq n-1~.
\la{lkskt0}
\eea
An alternative way to write this  condition is
\bea
\omega^\ell\wedge d(\omega^{k-1}\wedge H)=0~, \la{lkskt}
\eea
where the 3-form $H=-i (\partial-\bar\partial)\omega$  is  the
torsion of  $\hat\nabla$, and $\hat\nabla$   is the unique
Hermitian connection with skew-symmetric torsion.  There are some
advantages of writing (\ref{lkskt0}) as (\ref{lkskt})  as the
latter can be generalized to all G-structures which admit a
compatible connection with skew-symmetric torsion. These include
for example $U(n)$, $SU(n)$, $Sp(n)$,  $Sp(n)\cdot Sp(1)$, $G_2$
and $Spin(7)$ structures \cite{iibhor}.

Special cases of (\ref{lkskt0}) and (\ref{lkskt}) conditions have appeared before   in the literature.
First take $\ell=0$, and so  (\ref{lkskt0}) and (\ref{lkskt}) can be rewritten as
\bea
\partial\bar\partial \omega^k=0~,~~~~~d(\omega^{k-1}\wedge H)=0~,~~~1\leq k\leq n-1~,
\la{kskt} \eea respectively. For $k=n-1$, the above conditions
coincide with the {\it Gauduchon structure} on Hermitian manifolds
\cite{gauduchon} which is usually written as $\delta \theta=0$,
where $\theta=\delta\omega\circ J$ is the Lee form of the
Hermitian form. As a consequence of the Gauduchon's theorem in the
every conformal class, there is a Hermitian metric which satisfies
the Gauduchon condition. Thus every Hermitian manifold admits a
Gauduchon structure. Furthermore for $k=n-2$, the (\ref{kskt})
condition has been called  {\it astheno-K\"ahler} \cite{yaujost}
and it has been studied in the context of harmonic maps  and in
connection  with the extension of Siu's rigidity theorem to
non-K\"ahler complex manifolds. Recently, examples of such
manifolds have been given in \cite{FT}.

Another special case that has been extensively investigated for
many years  is (\ref{kskt}) for $k=1$,  $n>2$, or equivalently
$dH=0$. This coincides with the {\it strong structure} on
Hermitian manifolds \cite{hkt}  and has found many applications in
both physics, see eg \cite{hull, howe, strominger, hethor, hethor1} and
geometry, see eg \cite{IP1, IP2, salamon, swann, yau, yau1, FT1, streets}. For example
in type II string theory,  $H$ is identified with the 3-form field
strength. This is  required by construction to satisfy $dH=0$.
Recently Streets and Tian \cite{streets} introduced a hermitian
Ricci flow under which the pluriclosed or equivalently strong KT structure is
preserved. Viewing (\ref{kskt}) as a generalization of the strong
condition on a Hermitian manifold, we shall refer to it uniformly
as {\it k-strong} condition and the associated Hermitian manifolds
as {\it k-strong K\"ahler with torsion} or k-SKT for short.
Similarly, we shall refer to (\ref{lkskt0}), or equivalently to
(\ref{lkskt}), as the  {\it $(\ell\vert k)$-strong} condition and
to the associated Hermitian manifolds as admitting  a {\it
$(\ell\vert k)$-strong K\"ahler  with torsion}  structure or
$(\ell\vert k)$-SKT for short.

More recently, (\ref{lkskt0}) for $\ell=n-k-1$, the  {\it
generalized k-Gauduchon condition}, has been used by Fu, Wang and Wu
\cite{FWW} to prove a generalization of the Gauduchon theorem.
Examples of manifolds which satisfy the 1-Gauduchon condition have
been given in \cite{FWW} and \cite{FU}.

Apart from the condition (\ref{kskt}) above, Hermitian  manifolds
can also be restricted by a holonomy condition. This is usually
expressed as the requirement that one of the Hermitian connections
has reduced holonomy $G\subset U(n)$. In many investigations, the
holonomy condition is imposed in addition to conditions like
(\ref{lkskt}) on the Hermitian form. In many applications, see eg
\cite{strominger, bismut, IP1, IP2, gold, grant,  yau, yau1, FIUV, hethor, hethor1,
iibhor}, the holonomy condition is imposed on  the Hermitian
connection with skew-symmetric torsion $\hat\nabla$. Because of
this, we say that a 2n-dimensional Hermitian manifold is {\it
$(\ell\vert k)$-strong Calabi-Yau with torsion}, or equivalently
$(\ell\vert k)$-SCYT,   iff it is $(\ell\vert k)$-SKT and
$$
{\mathrm{hol}}(\hat\nabla)\subseteq SU(n)~. $$
If $\ell=0$, we
simply refer to such manifolds as k-SCYT. It is clear from this
that the Ricci form $\hat\rho$ of $\hat\nabla$ on $(\ell\vert
k)$-SCYT manifolds must vanish, $\hat\rho=0$, and consequently,
the $\omega$-trace $\hat b$ of $\hat\rho$ is also zero, $\hat
b=0$. Such manifolds have appeared before in the literature. In
particular, it has been shown in \cite{iibhor} that the
supersymmetric IIB black hole horizons are 8-dimensional 2-SCYT
manifolds and some examples have been constructed.

The purpose of this paper is twofold. First, we shall prove some
vanishing results for the Dolbeault cohomology of $(\ell\vert
k)$-SKT and $(\ell\vert k)$-SCYT manifolds. Then we shall give
some examples of manifolds with k-SKT and k-SCYT structures. One
of our  main results is a vanishing theorem on the plurigenera, $$
p_m(J)=\mathrm {dim}\, H^0_{\bar \partial}(X,K^m)~,$$ which is the
dimension of the number of holomorphic sections of the m-th power
of the canonical bundle. In particular, one has the following.
\begin{thrm}
Let $M$ be a compact  2n-dimensional $(\ell\vert k)$-SKT manifold
satisfying the condition
\[
 \hat b+ {n-k-1\over 3(n-2)} ||H||^2+{2
(k-1)\over n-2}||\theta||^2>0~,~~~~n\not=2. \]
Then
$$p_m(J)=0, m>0.$$
\end{thrm}
The proof of this result is based on an inequality derived in \cite{IP2} for the vanishing of plurigenera
for KT manifolds.

A consequence of this is the following.
\begin{thrm}
Let $M$ be a compact 2n-dimensional $(\ell\vert k)$-SCYT non Calabi-Yau manifold
and  $k<n-1$, $n>2$. Then  the plurigenera $p_m(J)=0, m>0$, and so
$M$ does not admit a holomorphic (n,0)-form.

In particular,
$p_m(J)=0, m>0$ for all generalized k-Gauduchon manifolds,
$k<n-1$, for which ${\rm hol}(\hat\nabla)\subseteq SU(n)$.
\end{thrm}
This theorem generalizes the results obtained in \cite{IP1, IP2,
FU} for 1-SCYT and $(n-2\vert 1)$-SCYT manifolds  and that
obtained in \cite{FT} for $(n-2)$-SCYT manifolds.

An immediate application of theorem 1.2 is that every conformally
balanced, $\theta=d\phi$, $\phi$ a function on $M$, compact
$(\ell\vert k)$-SCYT manifold, $k<n-1$, $n>2$, is Calabi-Yau. This
is similar to the result originally proved for the special case of
conformally balanced 1-SCYT and $(n-2\vert 1)$-SCYT manifolds
manifolds in \cite{IP1, IP2, FU} and adapted to $(n-2)$-SCYT
manifolds in \cite{FT}. However, one can generalize these results
using the   work of \cite{FWW} on generalized k-Gauduchon manifolds. In
particular, one has the following.
\begin{thrm}
Every compact, conformally balanced,  $(\ell\vert k)$-SKT
manifold, $k<n-1$, $n>2$, is  K\"ahler.
\end{thrm}

One can also consider  Hermitian manifolds admitting
a {\it generalized $(\ell\vert k)$-SKT structure} given by
\bea
2k i\, \omega^\ell\wedge
\partial\bar\partial \omega^k\equiv d(\omega^k\wedge H)\wedge \omega^\ell=
{1\over (k+\ell+1)!} \,\alpha_{k,\ell}\, \omega^{k+\ell+1}~, ~~~1\leq k+\ell\leq n-1~,
\la{glkskt0}
\eea
where $\alpha_{k,\ell}$ is a function on $M^{2n}$ which depends on $\omega$. Our results in theorems (1.1) and (1.2) generalize to this case
 provided
that $\alpha_{k,\ell}>0$.

Some of our results also apply to $(k_1,k_2,k_3)$-strong hyper-K\"ahler manifolds with torsion ($(k_1,k_2,k_3)$-SHKT)  \cite{iibhor}.
In particular, one can show that for all   $(\ell,1,1\vert k,1,1)$-SHKT manifolds with $k<n-1$ and $n>2$,  $p_m(I)=0, m>0$,
and cyclicly for $J$ and $K$. Furthermore if a $(\ell,1,1\vert k,1,1)$-SHKT manifold $M^{4n}$ is conformally balanced,
$\theta_{\omega_I}=2 d\phi$, $\phi$ a function on $M^{4n}$, then $M^{4n}$ is hyper-K\"ahler. The latter also applies
cyclicly for $J$ and $K$. These statements follow  because these $(k_1,k_2,k_3)$-SHKT structures
are special cases of the $(\ell\vert k)$-SCYT structures that appear in  theorems (1.1) and (1.2). Because of this,
we shall not elaborate further.

We shall construct several examples of k-SKT. Some examples of 2-SKT and 2-SCYT manifolds have already been given in \cite{iibhor}.
Here we shall extend a method initially used  by Swann  \cite{swann} to construct examples of
 1-SKT and HKT manifolds to give new examples of  k-SKT and k-SCYT manifolds. In particular, we shall construct several simply connected examples.

\newsection{Vanishing theorems for  2-SKT and 2-SCYT manifolds}

It is instructive to first prove theorems (1.1) and (1.2) for  2-SKT and 2-SCYT manifolds
and then extend the results to the most general case. In particular, this will
establish the results of theorems (1.1) and (1.2) for 2-Gauduchon manifolds.
As it has been mentioned in the introduction the theorems (1.1) and (1.2) have already been proven
for 1-SKT and 1-SCYT manifolds, respectively \cite{IP1, IP2}. Essentially the proof extends to 1-Gauduchon
manifolds, see also \cite{FU}. We shall demonstrate the proof of
theorem (1.3) after those for the theorems (1.1) and (1.2) for $(\ell\vert k)$-SKT and $(\ell\vert k)$-SCYT manifolds.  Before we proceed with this, we establish our conventions.

\subsection{Conventions and preliminaries}

Let  $(M,J,g)$   a Hermitian manifold of dimension 2n. Then, the Hermitian form\footnote{There is a sign difference
 from the definition of $\omega$ given in \cite{FWW} which is important in the proof of theorem (1.3).}
 is defined
as $\omega(X,Y):=g(X,JY)$ or equivalently in components
\[
\omega_{ij}= g_{ik} J^k{}_j~. \] The torsion $H$ of $\hat\nabla$
is $H=-i(\partial-\bar\partial)\omega$ or equivalently
\[
    H(X,Y,Z)=d^c\omega(X,Y,Z)=Jd\omega(X,Y,Z)=-d\omega(JX,JY,JZ)~,
 \]
 where we have use  that $JF(X_1,\dots,X_r):=(-1)^rF(JX_1,\dots,JX_r)$ for a r-form $F$.

 For the curvature we use the convention $\hat R=[\hat\nabla,\hat\nabla] -
 \hat\nabla_{[,]}$ Consequently, $\hat\rho(X,Y)=\hat
 R(X,Y,e_i,Je_i), \quad b=\hat\rho(Je_i,e_i)$, where we use Einstein summation conventions, ie repeated indices are summed over.

 The Lee form $\theta:=\delta\omega\circ J$ of the Hermitian manifold is given in terms of $H$ as
 \[
\theta(X)=-\frac12H(JX,e_i,Je_i)=\frac12g(H(JX),\omega)=\frac12(\omega\lrcorner
H(JX))~,~~~\theta_i={1\over 2} J^k{}_i H_{kj\ell} \omega^{j\ell}~.
\]
Moreover, we define the (1,1) form \cite{IP2}
\[
\lambda(X,Y):=dH(X,Y,e_i,Je_i)=-g(dH(X,Y),\omega)=-(\omega\lrcorner
dH(X,Y))~,~~
\]
i.e. $\lambda_{ij}=- dH_{ijk\ell} \omega^{k\ell}$.
We also write
\[ ||H||^2=H(e_i,e_j,e_k)H(e_i,e_j,e_k)=H_{ijk}
H^{ijk}~. \]
As a volume form, we use $d{\rm vol}(M)= {1\over
n!}\omega^n$, where $\omega^p= \wedge^p \omega$. In particular
${1\over p!}\star\omega^p={1\over (n-p)!}\, \omega^{n-p}$, where
$\star$ is the Hodge star operator.

\subsection{The $\alpha_2$ function}

On a 2n dimensional hermitian manifolds $(M,g,J)$, we define the function $\alpha_k(\omega)$ by
\bea
\alpha_k(\omega)=2k \star(i\partial\bar\partial\omega^k\wedge\omega^{n-k-1})= \star(d(\omega^{k-1}\wedge H)\wedge\omega^{n-k-1})~.
\la{g00}
\eea
Clearly, $\a_k(\omega)=0$  provided  $(M,g,J)$ admits a $(\ell\vert k)$-SKT structure.

Using $\omega^{n-k-1}={(n-k-1)!\over ( k+1)!}\star\omega^{k+1}$,  we  have the expression
\bea
\alpha_k(\omega)={(n-k-1)!\over ( k+1)!}\star(d(\omega^{k-1}\wedge H)\wedge\star\omega^{k+1})=
g(d(\omega^{k-1}\wedge H),\omega^{k+1})\nonumber\\=
(-1)^{k+1}{ (n-k-1)!\over (k+1)!\cdot 2^{k+1}}\,d(\omega^{k-1}\wedge H)(e_{i_1},Je_{i_1},\dots,e_{i_{k+1}},Je_{i_{k+1}})
\la{g0}
\eea

First we calculate $\a_2$. For this let us consider the following.
\begin{lemma}
On a 2n-dimensional hermitian manifold we have
\begin{equation}\label{01}
-\omega\lrcorner d(\omega\wedge H)=(8-2n)dH+\omega\wedge\lambda+2(JH\wedge J\theta)+2(H\wedge\theta)+2J(e_i\lrcorner H)\wedge(e_i\lrcorner H).
\end{equation}
In particular, on a 2-SKT manifold of dimension 2n we have:
\begin{equation}\label{1}
(8-2n)dH+\omega\wedge\lambda+2(JH\wedge J\theta)+2(H\wedge\theta)+2J(e_i\lrcorner H)\wedge(e_i\lrcorner H)=0.
\end{equation}
\end{lemma}
\emph{Proof.}
The identity $d(\omega\wedge H)=d\omega\wedge H+\omega\wedge dH$ also reads
\begin{equation}\label{2}
d(\omega\wedge H)=-JH\wedge H+\omega\wedge dH~.
\end{equation}
A straightforward calculation using our conventions reveals that
\bea
(-JH\wedge H)(X,Y,Z,U,e_i,Je_i)&=&[2JH\wedge J\theta+2H\wedge\theta\cr
&&~~~~~+[2J(e_i\lrcorner H)\wedge(e_i\lrcorner H)](X,Y,Z,U);\la{new1}
\\
(\omega\wedge dH)(X,Y,Z,U,e_i,Je_i)&=&[\omega\wedge\lambda
+(8-2n)dH](X,Y,Z,U). \la{new2} \eea
The last two equalities
together with \eqref{2} imply \eqref{01}.

The 2-SKT condition $d(\omega\wedge H)=0$ and \eqref{01} give \eqref{1}.
\hfill \emph{Q.E.D.}
\begin{cor} On a 2n-dimensional Hermitian manifold, we have
\bea\nonumber
\omega\lrcorner(\omega\lrcorner
d(\omega\wedge H))&=&(12-4n)\lambda+\lambda(e_i,Je_i)\omega
+8\theta\wedge J\theta \cr
 &&-8(J\theta)\lrcorner H-8J((J\theta\lrcorner H)
+4J(e_ie_j\lrcorner H)\wedge (e_ie_j\lrcorner H).
\eea
In particular, on a 2-SKT manifold of dimension 2n we have:
\begin{equation}\label{3}
(4n-12)\lambda=\lambda(e_i,Je_i)\omega +8\theta\wedge J\theta -8(J\theta)\lrcorner H-8J((J\theta\lrcorner H)
+4J(e_ie_j\lrcorner H)\wedge (e_ie_j\lrcorner H).
\end{equation}
\end{cor}
\emph{Proof.} Taking the traces in \eqref{01}, we get
\bea\nonumber
&&d(\omega\wedge
H)(e_i,Je_i,e_j,Je_j,X,Y)=(12-4n)\lambda(X,Y)+
\lambda(e_i,Je_i)\omega(X,Y) \cr &&~~~
-8\theta(X)\theta(JY)+8\theta(JX)\theta(Y)
-8H(X,Y,J\theta)-8H(JX,JY,J\theta) \cr
&&~~~-4H(JX,e_i,e_j)H(Y,e_i,e_j)+4H(X,e_i,e_j)H(JY,e_i,e_j)~, \eea
which proves the assertion. \hfill \emph{Q.E.D.}
\begin{prop}
On a 2n-dimensional Hermitian manifold the function $\a_2$ is given by
\bea
\a_2(\omega)=(n-3)!\Big[(n-2)\delta\theta+ (n-3)\big[||\theta||^2-{1\over6}||H||^2\Big]
\la{g1}
\eea
In particular, on a 2-nd Gauduchon manifold as well as on a 2-SKT manifold we have
\bea
(n-2)\delta\theta+(n-3)\Big[||\theta||^2-{1\over6}||H||^2\Big]=0.
\la{g3}
\eea
\end{prop}
\emph{Proof.}
The trace in \eqref{3} together with \eqref{g0}  gives
\begin{lemma}
On a 2n-dimensional hermitian manifold the function $\a_2(\omega)$  is given by
\bea
\a_2(\omega)={(n-3)!\over 2^3} \Big[(n-2)\lambda(e_i,Je_i)-8||\theta||^2+{4\over3}||H||^2\Big]
\la{g2}
\eea
On a 2-SKT manifold of dimension 2n we have:
\begin{equation}\label{5n}
(n-2)\lambda(e_i,Je_i)=8||\theta||^2-{4\over3} ||H||^2=0.
\end{equation}
\end{lemma}
To complete the proof of the proposition, we use the identity
\begin{equation}\label{6}
\lambda(e_i,Je_i)=8\delta\theta+8||\theta||^2-\frac43||H||^2~,
\end{equation}
established in \cite{AI,IP1} in the context of KT manifolds.  Combining \eqref{6} with \eqref{g2}, it is straightforward to prove  \eqref{g1}.\hfill \emph{Q.E.D.}

\begin{cor}
On a 2n-dimensional  2-SKT manifold, one has
\begin{equation}\label{7n}
(n-3)\lambda(e_i,Je_i)=-8\delta\theta~,
\end{equation}
and
\begin{equation}\label{7nn}
(n-2)\lambda(e_i,Je_i)=8||\theta||^2-{4\over3}||H||^2~.
\end{equation}
\end{cor}
\emph{Proof.} The proof of the above two equations follows from \eqref{5n} and \eqref{6}. \hfill \emph{Q.E.D.}

\subsection{Proof of theorems (1.1) and (1.2)}

{\it Proof of theorem (1.1):} Now let us turn to the proof of
theorem (1.1) for 2-SKT manifolds. For this, we use the result in
\cite[Theorem 4.1]{IP2} that the plurigenera, $p_m(J), m>0$, of a
KT manifold vanish provided that \bea \hat b+||C||^2+ \frac14
\sum^{2n}_{i=1}\lambda(e_i,Je_i)>0~, \label{ineq} \eea where $\hat
b$ is the $\omega$-trace of the Ricci form $\hat\rho$ of
$\hat\nabla$ and $C$ is the torsion of the Chern connection. The
fact that $H$ is of type (1,2)+(2,1) implies \bea
H(Je_k,Je_i,e_j)H(e_k,e_i,e_j)=\frac13||H||^2. \la{21} \eea We
recall that the torsion $C$ of the Chern connection of a KT
manifold $(M,g,J)$ is expressed in terms of  $H$  as, \[
g(C(X,Y),Z)=\frac12H(X,JY,JZ)+\frac12H(JX,Y,JZ)~,
 \]
 see e.g. \cite{IP1}.
Using this, \eqref{21} and   that  $H$ is a (1,2)+(2,1)-form, one finds that
\begin{equation}\label{chern}
||C||^2=\frac13||H||^2.
\end{equation}
Next using \eqref{chern} and \eqref{7nn}, one has that
\begin{eqnarray}\label{8n}
&&\hat b+||C||^2+\frac14\lambda(e_i,Je_i)=\hat b+\frac13 ||H||^2+{2\over n-2} ||\theta||^2-{1\over 3(n-2)}||H||^2
\cr
&&=\hat b+\frac13 (1-{1\over n-2})||H||^2+{2\over n-2}||\theta||^2>0
\end{eqnarray}
which is positive for $n>2$ according to the condition of theorem (1.1). This establishes theorem (1.1) for $k=2$. \hfill \emph{Q.E.D.}

\vskip 0.3cm {\it Proof of theorem (1.2):} Now, let us turn to
theorem (1.2) for 2-SCYT manifolds. It readily follows from
theorem (1.1). Since the holonomy of the  connection  with
skew-symmetric torsion $\hat\nabla$ is in $SU(n)$, $\hat b=0$, and the
inequality \eqref{8n} is always satisfied provided that  $H$ does
not vanish. Clearly the above statement also holds under the
weaker assumption that $\hat b=0$.\hfill \emph{Q.E.D.}

\begin{cor}
A compact, conformally balanced, 2-SCYT manifold is Calabi-Yau.
\end{cor}
{\it Proof :}  This is a special case of theorem (1.3) which we
shall demonstrate later.  This is also an extension of a similar
theorem proved in \cite{IP1} for conformally balanced 1-SCYT
manifolds. It follows from \cite{strominger}  that a
2n-dimensional  conformally balanced CYT manifold admits a
holomorphic (n,0)-form. Combining this with the statement of
theorem (1.2) for 2-SCYT manifolds, one concludes that $H=0$, and
so $M$ is Calabi-Yau.\hfill \emph{Q.E.D.}

\newsection{Vanishing theorems for 2n-dimensional k-SKT and  k-SCYT manifolds}

\subsection{The $\alpha_k$ function}

We have shown theorems (1.1) and (1.2)  for k-SKT and k-SCYT manifolds for $k=1,2$. It remains to extend
these to all $(\ell\vert k)$-SKT and $(\ell\vert k)$-SCYT manifolds for $k>2$, $\ell>0$ . Instrumental in this is the generalization of \eqref{g1} and \eqref{g3} for $k>2$.

\begin{prop}
On a 2n-dimensional Hermitian manifold the function $\a_k$ is given by
\bea
\a_k(\omega)=(n-3)!\Big[(n-2)\delta\theta+(n-k-1)\big[||\theta||^2-{1\over6}||H||^2\big]\Big]~.
\la{gk1}
\eea
In particular, a 2n-dimensional Hermitian manifold is  generalized k-Gauduchon, if and only if, the next identity holds
\bea
(n-2)\delta\theta+(n-k-1)\Big[||\theta||^2-{1\over 6}||H||^2\Big]=0.
\la{gk3}
\eea
\end{prop}
\emph{Proof.} First we show
\begin{lemma}
Let $M^{2n}$ be a  Hermitian  manifold. Then
\begin{equation}\label{5nk}
\a_k(\omega)={(n-3)!\over 2^3} \big [
(n-2)\sum_{i=1}^{2n}\lambda(e_i,Je_i)-8 (k-1) ||\theta||^2+{4\over3} (k-1)||H||^2\big]
\end{equation}

\end{lemma}
To prove the lemma, we write \eqref{g00} as
\bea
\alpha_k(\omega)&=&\star(d(\omega^{k-1}\wedge H)\wedge\omega^{n-k-1})=
\star\Big(\Big[(k-1) \omega^{k-2}\wedge d\omega\wedge H+ \omega^{k-1}\wedge dH\Big]\wedge\omega^{n-k-1}\Big)\nonumber
\\
&=&\star\Big(\omega^{n-3}\wedge \Big[(k-1) d\omega\wedge H+\omega\wedge dH\Big]\Big)=
\star\Big(\omega^{n-3}\wedge \Big[-(k-1)JH\wedge H+\omega\wedge dH\Big]\Big)\nonumber
\\
&=& {(n-3)!\over 3!}
\star\Big(\star\omega^3\wedge \Big[-(k-1)JH\wedge H+\omega\wedge dH\Big]\Big)
\nonumber
\\
&=&g(\omega^3,\Big[-(k-1)JH\wedge H+\omega\wedge dH\Big])\nonumber
\\
&=&-{(n-3)!\over 3!\cdot 2^3}\Big(-(k-1)JH\wedge H+\omega\wedge dH\Big)(e_i,Je_i,e_j,Je_j,e_k,Je_k)\nonumber
\\
&=& {(n-3)!\over 2^3} \big [
(n-2)\sum_{i=1}^{2n}\lambda(e_i,Je_i)-8 (k-1) ||\theta||^2+{4\over3} (k-1)||H||^2\big],
\la{g01}
\eea
where we used the \eqref{new1} and \eqref{new2} and their traces. This completes the proof of the lemma.

Next to prove the proposition,  substitute \eqref{6} into \eqref{5nk} to get \eqref{gk1}. This  completes the proof. \hfill \emph{Q.E.D.}

A generalization of \eqref{7n} and \eqref{7nn} is as follows.
\begin{cor}
Let $M$ be  a 2n-dimensional  $(\ell\vert k)$-SKT manifold, then
\begin{equation}\label{7nk}
(n-k-1) \sum_{i=1}^{2n}\lambda(e_i,Je_i)=-8(k-1) \delta\theta
\end{equation}
and
\begin{equation}\label{7nnk}
(n-2) \sum_{i=1}^{2n}\lambda(e_i,Je_i)=8(k-1) ||\theta||^2-{4\over3} (k-1) ||H||^2~.
\end{equation}
\end{cor}
\emph{Proof.} It follows immediately as an application of \eqref{6}, \eqref{5nk} and the fact that $\a_k(\omega)=0$ for
all $(\ell\vert k)$-SKT manifolds. \hfill \emph{Q.E.D.}

Integrate \eqref{gk3} over a compact $M$ observing that $||\theta||^2=||J\delta\omega||^2=||\delta\omega||^2$  and $||H||^2=||d\omega||^2$ to obtain
\begin{cor}
Let $(M,\omega)$ be  a compact 2n-dimensional  $(\ell\vert k)$-SKT manifold, then for $k<1<n-1$ we have
\[
\int_M||\delta\omega||^2dvol(M)=\int_M\frac16||d\omega||^2dvol(M).
\]
\end{cor}

\subsection{Proof of Theorems (1.1) and (1.2)}

 {\it Proof of theorem (1.1):} To show this  for all $(\ell\vert k)$-SKT manifolds, we apply again the inequality \eqref{ineq} established
in \cite[Theorem 4.1]{IP2} as a condition for the vanishing of plurigenera for KT manifolds and use \eqref{7nnk}. One finds that

\bea\label{8nk}
&&\hat b+||C||^2+\frac14\lambda(e_i,Je_i)=\hat b+\frac13 ||H||^2+{2(k-1)\over n-2} ||\theta||^2-{k-1\over 3(n-2)}||H||^2
\cr
&&=\hat b+{n-k-1\over 3(n-2)} ||H||^2+{2 (k-1)\over n-2}||\theta||^2>0~.
\eea
which is positive for $n>2$ according to the condition of theorem (1.1).
\hfill \emph{Q.E.D.}

\vskip 0.3cm

{\it Proof of theorem (1.2):} Now if $M$ is  $(\ell\vert k)$-SCYT,
then one has that $\hat b=0$. This follows from the requirement that
the holonomy of the Hermitian  connection with skew-symmetric
torsion, $\hat\nabla$, is contained in $SU(n)$. It is clear then
that the inequality \eqref{8nk} always holds and so $p_m(J)=0,
m>0$ for all  $(\ell\vert k)$-SCYT manifolds.\hfill \emph{Q.E.D.}


Theorems (1.1) and (1.2) can be extended to the generalized $(\ell\vert k)$-SKT and $(\ell\vert k)$-SCYT manifolds as follows.
\begin{cor}
Let $M^{2n}$ be a non K\"ahler generalized $(\ell\vert k)$-SKT
manifold, then $p_m(J)=0$, $m>0$, provided that
\bea\nonumber \hat
b+ {n-k-1\over 3( n-2)}||H||^2+{2 (k-1)\over n-2}||\theta||^2+ {2
n (n-1)\over (k+\ell+1)!}  \alpha_{k,\ell}>0~,~~~~n\not=2, \eea
where $\alpha_{k,\ell}$ is given in \eqref{glkskt0}. In
particular, the plurigenera vanish for every generalized
$(\ell\vert k)$-SCYT manifold for which $\alpha_{k,\ell}\geq 0$.
\end{cor}
The proof of this follows immediately from those of theorems (1.1) and (1.2) above. Note that
$\alpha_k={n!\over (k+\ell+1)!} \alpha_{k,\ell}$.

\begin{cor}
A compact, conformally balanced, $(\ell\vert k)$-SCYT manifold is
Calabi-Yau.
\end{cor}
{\it Proof :} This is a special case of theorem (1.3) and it
follows directly from the results of \cite{IP1} together with
theorems (1.1) and (1.2).  The proof is similar to that given  as
for the case of  conformally  balanced 2-SCYT manifolds in section
2.\hfill \emph{Q.E.D.}

\subsection{Proof of Theorem (1.3)}

It has been shown in \cite{FWW} that on a compact Hermitian
manifold there is a unique constant $\gamma_k(\omega)$ invariant
under biholomorphisms which depends smoothly on $\omega$ such that
the k-generalized Gauduchon equation\footnote{The sign difference
in \eqref{kgc} from that in \cite{FWW} is conventional and it is
due to a sign difference in the definition of Hermitian form
$\omega$.} \bea
 {i\over2}  e^{-u}\partial\bar\partial (e^u
\omega^k)\wedge \omega^{n-k-1}=-\gamma_k(\omega)\, \omega^n~,
\la{kgc} \eea has  a solution $u$, where $u$ is uniquely
determined up to a constant. In particular, a Hermitian manifold
$M$ admits a generalized k-Gauduchon metric in the conformal class
of $\omega$,  if and only if $\gamma_k=0$ \cite[Proposition
8]{FWW}.

The existence of generalized k-Gauduchon metrics depends crucially on the sign
of $\gamma_k$. It is also shown \cite[Proposition 11]{FWW} that the sign of $\gamma_k(\omega)$
remains constant in the conformal class of $\omega$. Moreover,
\cite[Proposition 12]{FWW}, in our notations, tells us that
$\gamma_k(\omega)>0 (=0, or  <0)$ if there exists a hermitian form
$\omega'$ in the conformal class of $\omega$ such that
$\alpha_k(\omega')<0 (=0, or
>0),$ respectively.

Suppose now that $\omega$ is conformally balanced. In such case,
there is a function $f$ on $M$ specified up to a constant such
that $\tilde\omega=e^{f} \omega$ is balanced, ie the corresponding
Lee form $\tilde\theta=0$. The next lemma makes \cite[Lemma
16]{FWW} more precise and proofs our Theorem 1.3.
\begin{lemma} On a compact balanced non-K\"ahler Hermitian manifold $(X,\tilde\omega)$ the
constant $$\gamma_k(\tilde\omega)>0, \quad for \quad 1\le k\le
n-2$$
\end{lemma}
\emph{Proof.} To proof the lemma  substitute $\tilde\theta=0$ into
\eqref{gk1} to conclude
$$\alpha_k(\tilde\omega)=-(n-3)!\frac{n-k-1}6||H(\tilde\omega)||^2.$$
Therefore for $k\not=n-1$, $\alpha_k(\tilde\omega)<0$, provided
that $H(\tilde\omega)\not=0$. Hence $\gamma_k(\tilde\omega)>0$ and the lemma follows. \hfill
\emph{Q.E.D.}

To complete the proof of the theorem 1.3 recall that the sign of
$\gamma_k$ does not depend on the conformal class of $\omega$ and
if $H(\tilde\omega)\not=0$, then also $\gamma_k(\omega)>0$.
 Now from the assumptions of theorem (1.3), $\omega$ is
$(\ell\vert k)$-SKT and therefore generalized k-Gauduchon which requires that
$\gamma_k(\omega)=0$. This leads to a contradiction unless
$H(\tilde \omega)=0$ and so $\tilde \omega$ is K\"ahler which
completes the proof of the theorem. \hfill \emph{Q.E.D.}

\subsection{Locally conformally K\"ahler manifolds} It is observed in \cite{AI} that 1-SKT
locally conformally K\"ahler manifold must be K\"ahler. Recently,
it is shown in \cite{FWW} that the standard hermitian structure on
$S^5\times S^1$ which is locally conformally K\"ahler has
$\gamma_1<0$. We show that this is true in general.

We recall that a Hermitian manifold $(X,\omega)$ is \emph{locally
conformally K\"ahler} if there locally exists a conformal metric
which is K\"ahler and this is not true globally. For $n>2$ this
condition is equivalent to the equation
$d\omega=\frac1{n-1}\theta\wedge\omega$ which, in terms of $H$,
reads
\begin{equation}\nonumber
H=\frac1{n-1}J\theta\wedge\omega
\end{equation}
We have
\begin{prop} On a compact locally conformally K\"ahler  2n-manifold $(X,\omega)$ the
constant $\gamma_k(\omega)$ is negative for $1\le k\le
n-2$,
$$\gamma_k(\omega)<0, \quad for \quad 1\le k\le
n-2.$$
In particular, compact locally conformally K\"ahler structure does not admit $(l|k)$-SKT structure.
\end{prop}
\emph{Proof.} Let $\tilde\omega$ be the Gauduchon structure
globally conformal to $\omega$, which, in particular is locally
conformally K\"ahler and not K\"ahler. Then  we have
\begin{equation}\la{glck}
\tilde\delta\tilde\theta=0, \quad \tilde
H=\frac1{n-1}J\tilde\theta\wedge\tilde\omega,\quad
||\tilde\theta||^2\not=0~,
\end{equation}
where $\tilde\theta$ and $\tilde H$ are the Lee form and 3-form torsion associated to $\tilde \omega$, respectively.
A straightforward calculation yields
\begin{equation}\la{glck1}
||\tilde H||^2=\frac6{n-1}||\tilde\theta||^2.
\end{equation}
To proof the assertion  substitute \eqref{glck} and \eqref{glck1} into
\eqref{gk1} to conclude
$$\alpha_k(\tilde\omega)=(n-3)!(n-k-1)\Big[||\tilde\theta||^2-\frac16||H(\tilde\omega)||^2\Big]=
(n-3)!(n-k-1)\frac{n}{n-1}||\tilde\theta||^2.$$
Therefore for $k\not=n-1$, $\a_k(\tilde\omega)>0$. Consequently, $\gamma_k(\tilde\omega)<0$. \hfill
\emph{Q.E.D.}

\newsection{Fibrations and k-SKT structures}
\subsection{ k-SKT structures on product manifolds}

We shall focus on the construction of k-SKT and k-SCYT structures are they are more restrictive than $(\ell\vert k)$-SKT
and $(\ell\vert k)$-SCYT, respectively. In particular, if a Hermitian manifold admits a k-SKT or k-SCYT structure, then it also
admits a $(\ell\vert k)$-SKT or $(\ell\vert k)$-SCYT for all $\ell$.
It is straightforward to construct k-SKT structures on products of manifolds. In particular one has the following.
\begin{prop}
The product $M^{2m}\times N^4$ where $M^{2m}$  is a K\"ahler manifold and $N^4$ is Hermitian 4-manifold admits a k-SKT structure for all $k$.
\end{prop}
\emph{Proof.} Let $\omega_{(4)}$ be the Hermitian form of a Gauduchon structure on $N^4$. Then $N$ is an 1-SKT manifold
with respect to $\omega_{(4)}$, ie $dH_{(4)}=0$ as this coincides with the co-closure of the Lee form.
If  $\omega_{(2m)}$ is the K\"ahler form on  $M^{2m}$, then
\[
d((\omega_{(2m)}+\omega_{(4)})^k\wedge H_{(4)})= d(\omega_{(2m)}^k
\wedge H_{(4)}) =\omega_{(2m)}^k \wedge dH_{(4)}=0~. \]
This
proves the proposition. \hfill  \emph{Q.E.D.}

For an explicit example one can take $N^4=S^1\times S^3$ and $M^{2m}=\CP^m$.

Similarly, it is straightforward to see the following.
\begin{prop}
Let $N$ be a k-SKT manifold for $k\leq \ell$. Then the product $M^{2m}\times N$, where $M^{2m}$  is a K\"ahler manifold, is also k-SKT manifold  for all $k\leq\ell$.
\end{prop}

\subsection{The Swann twist}

KT and CYT manifolds can be constructed using torus fibrations, see \cite{sethi, gold, grant}. These provide a large
class of examples and so some of them may admit the more restrictive k-SKT and k-SCYT structures. Although this can be done directly by consider torus fibrations over suitable base spaces, it is advantageous  to use
 a  construction  proposed by Swann \cite{swann} to find 1-SKT and (1,1,1)-SHKT metrics. This will be adapted  to give new examples  k-SKT and k-SCYT manifolds. We begin with a summary of the Swann's twist construction.

 Let $M^{2n}$ be a Hermitian manifold equipped with a $T^p$ torus action $A_M$ which preserves the Hermitian structure.
 Denote the Lie algebra of the group $T^p$ acting on $M^{2n}$ with $\mathfrak{a}_M$.
 In addition, let $P$ be a $T^p$ principal bundle over $M^{2n}$ equipped with a  connection $\lambda$. Clearly
 $\lambda\in \Omega^1(P, {\mathfrak{a}}_P)$, where ${\mathfrak{a}}_P$ is the Lie algebra of $T^p$ which acts on $P$ from the right. Suppose
 now that the $A_M$ group action on $M^{2n}$ can be lifted to and a $T^p$ action $A_P$ on $P$. If $\xi$'s are the vector fields
 generated by the action of $A_M$  on $M^{2n}$, then the $A_P$ action on $P$ is generated by the vector fields
 \bea\nonumber
 \mathring \xi=\tilde \xi+ \mathring{\beta}\,\rho~,
 \eea
 where $\tilde \xi$ is the horizontal lift of $\xi$ with respect to $\lambda$, ie $\lambda(\tilde\xi)=0$, $\rho$ are the vectors
  generated by the right action of $T^p$ on $P$ and $\beta\in \Omega^0(P, \mathfrak{a}_P\otimes \mathfrak{a}^*_M)$. Necessary
  conditions for  the $T^p$ action on $M^{2n}$ to lift to $P$ are
  \bea
  {\cal L}_\xi F=0~,~~~i_\xi F=d\beta,~~~i_\xi i_\xi F=0~,
  \la{liftcon}
  \eea
  and $\mathring{\beta}=\pi^* \beta$, where $\beta\in \Omega^0(M, \mathfrak{a}_P\otimes \mathfrak{a}^*_M)$ with ${\cal L}_\xi \beta=0$,  $\pi$ is the projection of $P$ onto $M^{2n}$ and $\pi^* F=d\lambda$ is the curvature
  of $\lambda$.

  Provided that (\ref{liftcon}) holds, there is a lift of the $A_M$ action to $P$ which covers $A_M$ and commutes
  with the right action on $P$. This lift is not unique
  because $\beta$ is determined up to a constant $\nu$. For every choice $\nu\in \Omega^0(M, \mathfrak{a}_P\otimes \mathfrak{a}^*_M)\otimes\bZ$, one finds another lift of the $A_M$ action. All these lifts are free provided that
  $A_M$ group action on $M^{2n}$ is free.

 The Swann twist is  a new fibration which is constructed by taking the quotient of $P$ with respect to $A_P$
 of  $T^p$. If $A_P$ is a free action, then $W=P/A_P$ is a manifold. Otherwise, it may have orbifold singularities.
 For the explicit examples we consider, $A_P$ is a free action. Under certain conditions, the Hermitian structure on $M^{2n}$
 can be inherited on $W^{2n}$. For this, one should induce a metric and a Hermitian form on $W^{2n}$ from those on $M^{2n}$. Let us first
 begin with forms. Given a form $\tau\in \Omega^\ell(M)$, one can define $\pi^*\tau\in \Omega^\ell(P)$. The aim is to construct
 a new form $\mathring \tau \in  \Omega^\ell(P)$ such that $\mathring \tau=\pi^*_W \tau_W$, where $\pi_W$ is the projection of $P$ onto $W^{2n}$. For this assume that $\beta$ is invertible and  take
 \bea
 \mathring \tau=\pi^*\tau- \lambda_{\beta^{-1}}^A \wedge \pi^*(i_{\xi_A} \tau)-...(-1)^{p(\ell)}{1\over \ell!}
 \lambda_{\beta^{-1}}^{A_1}\wedge \cdots\wedge \lambda_{\beta^{-1}}^{A_\ell} \pi^*(i_{\xi_{A_1}}\cdots i_{\xi_{A_\ell}} \tau)~,
 \la{flift}
 \eea
where  $\lambda_{\beta^{-1}}=\beta^{-1}\lambda\in \Omega^1(P, \mathfrak{a}_M)$, and  $p(\ell) =1$ if $[\ell/4]=1,2$ and $p(\ell) =0$ if $[\ell/4]=3,0$.
One  can verify that $i_{\mathring \xi} \mathring \tau=0$ and ${\cal L}_{\mathring \xi} \mathring \tau=0$
provided that ${\cal L}_{ \xi}\tau=0$. Therefore $\mathring \tau$ projects down onto $W^{2n}$, ie there is a $\tau_W$ such that
$\mathring \tau=\pi_W^* \tau_W$.

 Observe that to determine $\tau_W$ it suffices to know $\mathring \tau$ up to $\lambda$-terms. This is because
 all the components of $\mathring \tau$ proportional to $\lambda$'s are determined by the $\lambda$ independent term
 and the vector fields $\xi$. In \cite{swann}, this is referred as $\mathscr{H}$-relation or equivalence, where $\mathscr{H}={\rm Ker}\,\lambda$. Because of this, it suffices to establish the various relations up to $\mathscr{H}$-equivalence. Suppose
 now that $\pi^*_M \tau =_{\mathscr{H}} \chi$, then a direct application of (\ref{flift}) reveals that
 \bea
 \pi^*_M d\tau =_{\mathscr{H}} d\chi- {\cal F}^A \wedge \pi^*(i_{\xi_A} \tau)~,
 \la{dlift}
 \eea
where ${\cal F}=\beta^{-1} F\in \Omega^2(M, \mathfrak{a}_M)$.
 Using this construction, we can lift to $P$ both the metric and Hermitian form of $M^{2n}$ as
 \bea
 \mathring g&=& \pi^* g- 2\lambda_{\beta^{-1}}^A\otimes \pi^* \eta_A+ g(\xi_A, \xi_B) \lambda_{ \beta^{-1}}^A \otimes \lambda_{\beta^{-1} }^B ~,
 \cr
 \mathring \omega&=&\pi^*\omega-\lambda_{\beta^{-1}}^A \wedge \pi^*(i_{\xi_A} \omega)-{1\over2}\,
 \lambda_{\beta^{-1}}^A \wedge \lambda_{\beta^{-1}}^B\,\, \pi^*(i_{\xi_A} i_{\xi_B}\omega)~,
 \la{wgo}
 \eea
  where $\eta_A(X)= g(\xi_A, X)$.
 Provided that both $g$ and $\omega$ are invariant under $A_M$, these can be projected down to $W^{2n}$
 to define an almost Hermitian structure on $W^{2n}$.
 Using the $\mathscr{H}$-equivalence equivalence one can write
 \bea\nonumber
 \pi^*_W g_W=_{\mathscr{H}}\pi^* g~,~~~ \pi^*_W \omega_W=_{\mathscr{H}}\pi^* \omega~,~~~
 \eea

 It remains to find the conditions for the almost complex structure on $W^{2n}$ to be integrable. For this let $\mathscr{A}$ the set
 of all Killing vector fields in $M^{2n}$ and $\mathscr{A}_I=\mathscr{A}\cap I\mathscr{A}$. Clearly one can find basis $e_1,\dots e_{2s}, e_{2s+1},\dots e_k$ in $\mathscr{A}$
 which is an extension of the basis $e_1, \dots e_{2s}$ of  $\mathscr{A}_I$ with $I(e_{2j-1})=e_{2j}$. Next choose a
 basis $\epsilon^\a$ of (1,0)-forms in $M^{2n}$, lift them to $P$ and define $\mathring \epsilon^\a$. Then
 \bea\nonumber
 d \mathring\epsilon^\a= \pi^* d \epsilon^\a- {\cal  F}\, \pi^* i_\xi \epsilon^\a~.
 \eea
 The complex structure on $W^{2n}$ is integrable iff the (0,2)-part of the above 2-form vanishes. The (0,2) component of
 $d \epsilon^\a$ vanishes as consequence of the integrability of the complex structure on $M^{2n}$. In addition the (0,2) component
 of the term involving $F$ also vanishes provided that
 \bea\nonumber
 ({\cal F}_{2j-1}+i {\cal F}_{2j})^{0,2}&=&0~,~~~~j=1,\dots, s~,
 \cr
{\cal  F}_r^{2,0}={\cal F}^{0,2}_r&=&0~,~~~~r=2s+1,\dots, k~,
 \eea
 where $F_i={\cal F}(e_i)$. In particular the complex structure on $W^{2n}$ is always integrable if $F$ is a (1,1)-form on $M^{2n}$.

 It remains to determine the torsion $H_W$. For this observe that
 \bea\nonumber
 \pi_W^* H_W=_{\mathscr{H}} \pi^*H+  i_I {\cal F}^A \wedge \pi^*(i_{\xi_A} \omega)-  {\cal F}^A \wedge \pi^*\eta_A~,
 \eea
 where again $\eta_A(X)= g(\xi_A, X)$. This follows directly from  (\ref{wgo}) and (\ref{dlift}) using $H=-i (\partial-\bar\partial)\omega$.
 In particular, if $F$ is a (1,1)-form, this simplifies to
 \bea\nonumber
 \pi_W^* H_W=_{\mathscr{H}} \pi^*H-
 {\cal  F}^A \wedge \pi^*\eta_A~.
 \eea

 \subsection{k-SKT structures  from K\"ahler manifolds}

 As a starting point let us take $X$ to be a K\"ahler manifold with K\"ahler form $\omega_X$ and take $M=X\times T^{2m}$. Assuming
 that $T^{2m}$ is also K\"ahler with respect to the standard flat metric and complex structure, clearly $M$ is
 a K\"ahler manifold with K\"ahler form $\omega_X+\omega_T$, where $\omega_T$ is the K\"ahler form of $T^{2m}$. Next assume that the $A_M$ action on $M$ is the standard one,  where the vector fields $\xi$
 generate the standard basis in the homology of $T^{2m}$. Choose now a principal $T^{2m}$ bundle over $X$ with $F$ a (1,1)-form on $X$. By construction $i_\xi F=0$, and so $\beta$ is a constant matrix. Furthermore, the condition $\beta\in \Omega^0(M, \mathfrak{a}_P\otimes \mathfrak{a}^*_M)\otimes\bZ$ implies that $\beta\in(\Lambda_P\otimes \Lambda^*_M)\otimes\bZ$, where
 $\Lambda_P$ and $\Lambda_M$ are the lattices used to construct the tori of the typical fibre of $P$ and that of the torus
 action on $M$.

 Next using some constant invertible matrix $\beta$, let us perform a Swann twist to find
 \bea\nonumber
 \pi_W^* H_W=_{\mathscr{H}}- {\cal  F}^A \wedge \pi^*\eta_A~.
 \eea
 The k-SKT condition on $W$ is satisfied provided\footnote{The inner derivation $i_L\chi$, of a vector k-form $L$ with a p-form $\chi$ is defined as
$
i_L\chi={1\over k!\cdot (p-1)!} L^j{}_{i_1\dots i_k} \chi_{ji_{k+1}\dots i_{p+k-1}} dx^{i_1}\wedge\cdots \wedge dx^{i_{p+k-1}}$.} that
 \bea
&&{\cal F}^A\wedge {\cal  F}^B \wedge i_{\xi^A}
 (\omega_X+\omega_T)^{k-1} \wedge \pi^*\eta_B
 \cr
 &&+{\cal  F}^A\wedge {\cal  F}^B\wedge (\omega_X+\omega_T)^{k-1}
 g_T(\xi_A, \xi_B)=0~.
 \la{ksktt}
 \eea
 For $k=1$, this becomes
 \bea\nonumber
{\cal  F}^A\wedge {\cal  F}^B\,
 g_T(\xi_A, \xi_B)=0~,
 \eea
which is the SKT condition derived in \cite{swann}.

 \subsubsection{2-SKT manifolds}

 Let us now consider the $k=2$ case. The condition (\ref{ksktt}) becomes
 \bea
  g_T(\xi_A, \xi_B)\,{\cal F}^A\wedge {\cal F}^B\wedge \omega_X=0~,
 \cr
 {\cal  F}^A\wedge {\cal  F}^B \wedge i_{\xi^A}
 \omega_T \wedge \pi^*\eta_B+{\cal  F}^A\wedge {\cal  F}^B\wedge \omega_T\,\,
 g_T(\xi_A, \xi_B)=0~.
 \la{2kskt}
 \eea
 \begin{prop}
 If $M=X\times T^2$ and $ g_T(\xi_A, \xi_B)=\delta_{AB}$, then $W$ is  2-SKT
 ifand only if
 \bea\nonumber
\delta_{AB}\,{\cal F}^A\wedge{\cal F}^B\wedge \omega_X =0~,~~~
 \eea
 \end{prop}
 {\it Proof }: It is easily seen that the second condition (\ref{2kskt}) is automatically satisfied. The first
 condition in (\ref{2kskt}) gives the restriction stated above.
  \hfill  \emph{Q.E.D.}

 It is therefore clear that if $W$ is SKT, then it is also 2-SKT. However, there are 2-SKT structures which are not
 induced from an 1-SKT one. To find one such example, take
 $X^6=X^4\times X^2$ and $m=2$, where $X^4$ and $X^2$ are K\"ahler manifolds with K\"ahler forms
 $\omega_{(4)}$ and $\omega_{(2)}$, respectively. Moreover, we choose $g_T(\xi_A, \xi_B)=\delta_{AB}$ as it is required in above proposition but take
\begin{equation}
 \beta^{-1}=\begin{pmatrix} p_1 & q_1\\ p_2 & q_2 \end{pmatrix}~,
 \la{sktalpha}
\end{equation}
  and $F^1$ to have support over $X^4$ while $F^2$ to have support over $X^2$ with $[F^1]\in H^2(X^4, \bZ)$ and $[F^2]\in H^2(X^2, \bZ)$.
  The 2-SKT condition then becomes
 \bea
( p_1^2+p_2^2) F^1\wedge F^1\wedge \omega_{(2)}+2 (p_1 q_1+p_2 q_2) F^2\wedge F^1\wedge \omega_{(4)}=0~.
\la{2sktf}
 \eea
 The condition that $\beta\in(\Lambda_P\otimes \Lambda^*_M)\otimes\bZ$ implies that $p_1, p_2, q_1, q_2$ are integers up
 to possibly multiplying them with $\det\beta$.

{\it Example~1:} There are many solutions to this equation. First,
suppose that $F^1$ can be chosen such that $F^1\wedge F^1=0$. Such
classes exist on any complex manifold $N$ which admits a
non-trivial holomorphic map $\Phi:~~N\rightarrow \CP^1$. Then
$F^1=\Phi^* \zeta$, where $\zeta\in H^2(\CP^1, \bZ)$. In
particular, $K_3$ admits two such representatives in the second
cohomology. For this one uses the Weierstrass $\wp$-function.
Other examples include  any 4-dimensional  K\"ahler manifold which
arises as a blow up at the intersection points of an algebraic
4-dimensional K\"ahler manifold with  a complex co-dimension $r-2$
hyperplane in $\CP^r$, see \cite{swann} for further explanation.
In such a case, the condition (\ref{2sktf}) reduces to
\bea\nonumber (p_1 q_1+p_2 q_2)=0~. \eea Using the scale
invariance of the equation, set $p_1=1$. Then $q_1=-p_2 q_2$ where
$p_2, q_2$ are any integers. The only additional requirement is
that $p_1 q_2-p_2 q_1\not=0$ for $\beta$ to be invertible. For an
explicit example set $X^4=K_3$ and $X^2=\CP^1$ with
$F^1=\wp^*\zeta$ and $F^2=\omega_{\CP^1}$. Then $W$ which is
identified as a $T^2$ bundle over $K_3\times \CP^1$ admits a 2-SKT
structure.

{\it Example~2:} For another example assume that
$[\omega_{(2)}]\in H^2(X^2, \bZ)$ and  $[\omega_{(4)}]\in H^2(X^4,
\bZ)$. Then set $F^1=\omega_{(4)}$ and $F^2=\omega_{(2)}$. The
resulting equation reads \bea\nonumber p_1^2+p_2^2+2 (p_1 q_1+p_2
q_2)=0~. \eea One solution to the above equation is $p_1=0$ and
$p_2+2 q_2=0$. Then $q_1\not= 0$ can be arbitrary. The only
additional condition is that $p_2\not=0$ which is required for
$\beta$ to be invertible.

Clearly there are many explicit examples by taking $X^2$ to be
$\CP^1$ and $X^4$ to a K\"ahler 4-dimensional manifold like
$\CP^2$. The resulting 8-dimensional 2-SKT manifold $W$ has at
most finite fundamental group.  So its universal cover $\tilde W$
will provide an example of a simply connected compact 2-SKT
manifold. In particular, $S^3\times S^5$ admits a 2-SKT structure.

\subsubsection{k-SKT manifolds}

The results described in the previous section can be generalized to k-SKT manifolds.
\begin{prop}
 If $M=X\times T^2$ and $ g_T(\xi_A, \xi_B)=\delta_{AB}$, then $W$ is  k-SKT iff
 \bea\nonumber
\delta_{AB}\,\,{\cal  F}^A\wedge {\cal  F}^B\wedge \omega^{k-1}_X=0~,~~~
 \eea
 \end{prop}
 The proof of this is similar to that given for 2-SKT manifolds.

 To find examples take $X=X^{2k}\times X^2$ with K\"ahler forms $\omega_{(2k)}$ and $\omega_{(2)}$, respectively,
 Take again a $T^2$ bundle over $X$ with curvature $(F^1, F^2)$ which has support on $X^{2k}$ and $X^2$, respectively.
 Then the 2-SKT condition reads
 \bea
 (k-1) (p_1^2+p_2^2) F^1\wedge F^1\wedge \omega_{2k}^{k-2}\wedge \omega_{(2)} +2 (p_1q_1+p_2q_2) F^1\wedge F^2\wedge \omega_{(2k)}^{k-1}=0~,
 \la{ksktf}
 \eea
 where we have chosen $\beta$ as in (\ref{sktalpha}).

{\it Example~1:} Clearly if $F^1\wedge F^1=0$, as in the example
given for the 2-SKT case in the previous section, the above
condition reduces to $(p_1q_1+p_2q_2)=0$. This is again solved as
in the 2-SKT case.

{\it Example~2:} Another possibility is to assume that
$[\omega_{(2)}]\in H^2(X^2, \bZ)$ and  $[\omega_{(2k)}]\in
H^2(X^{2k}, \bZ)$, and  set $F^1=\omega_{(2k)}$ and
$F^2=\omega_{(2)}$, then one finds that (\ref{ksktf}) reduces to
\bea (k-1) (p_1^2+p_2^2)+2 p_1 q_1+2p_2 q_2=0~.
\la{ksktex2}
\eea
A solution to this equation is $p_1=0$ and $(k-1) p_2+2 q_2=0$ for
arbitrary $q_1$. There are many solutions to these equations in
$\bZ$ for which $p_2\not=0$ which is required for $\beta$ to be
invertible. Taking $X^{2k}=\CP^k$ and $X^2=\CP^1$, one can show
that $S^{2k+1}\times S^3$ admits a k-SKT structure.

To give more examples take $M=X\times T^2$ with $ g_T(\xi_A,
\xi_B)=\delta_{AB}$ as before but now $X=X^{2k}\times X^4$, where
$X^{2k}$ and $X^4$ are K\"ahler manifolds with K\"ahler forms
$\omega_{(4)}$ and $\omega_{(2k)}$, respectively. Furthermore
assume that the $T^2$ fibration over $X\times T^2$ has curvature
$(F^1, F^2)$, where $F^1$ and $F^2$ have support on $X^{2k}$ and
$X^4$, respectively.

\begin{prop}
  $W$ admits a (k+1)-SKT condition provided
 \bea
{k (k-1)\over 2} (p_1^2+p_2^2) F^1\wedge F^1\wedge \omega_{(2k)}^{k-2}\wedge \omega_{(4)}^2&+&
2k (p_1 q_1+p_2 q_2) F^1\wedge F^2\wedge \omega_{(2k)}^{k-1}\wedge \omega_{(4)}
\cr
&+& (q_1^2 +q_2^2) F^2\wedge F^2\wedge\omega_{(2k)}^{k}=0~,
 \la{k1sktk}
 \eea
 where $\beta$ is chosen as in (\ref{sktalpha}).
 \end{prop}

{\it Example~3:} To find solutions to (\ref{k1sktk}) suppose that $F^1\wedge F^1=F^2\wedge F^2=0$. Then the condition reduces to requiring that
 $p_1 q_1+p_2 q_2=0$ which can be solved as in the 2-SKT case. For an explicit example take $M=K_3\times K_3\times T^2$
 and $F^1=\wp^*_1\zeta$ and $F^2=\wp^*_2\zeta$, where $\wp_1$ and $\wp_2$ are the Weierstrass functions of the first and second
 $K_3$ subspaces in $M$, respectively. This will give 3-SKT structures on  $T^2$ bundles over $K_3\times K_3$.

{\it Example~4:} Next assume that $[\omega_{(4)}]\in H^2(X^4,
\bZ)$ and  $[\omega_{(2k)}]\in H^2(X^{2k}, \bZ)$, and set
$F^1=\omega_{(2k)}$ and $F^2=\omega_{(2)}$.  Then substituting in
(\ref{k1sktk}),  one finds that $W$ admits a (k+1)-SKT condition
if \bea\nonumber {k (k-1)\over 2} (p_1^2+p_2^2) +2k (p_1 q_1+p_2
q_2)+ (q_1^2 +q_2^2)=0~. \eea

To find a solution set $p_1=0$ and observe that the above equation
can be rewritten as \bea\nonumber q_1^2+(k p_2+q_2)^2={k
(k+1)\over 2} p_2^2~. \eea This has solutions, eg $k=4$, $p_2=2$,
$q_1=2$ and $q_2=-2$.

\subsection{k-SCYT structures from K\"ahler-Einstein manifolds}

Examples of 2-SCYT manifolds have been constructed in
\cite{iibhor} in the context of IIB black hole horizons. Some of
the k-SKT manifolds we have constructed also admit a k-SCYT
structure. For this, we shall find the conditions such that the
Ricci form, $\hat\rho_W$, of the  connection $\hat\nabla_W$ with
skew-symmetric torsion on $W$  vanishes, $\hat\rho_W=0$. This
condition is equivalent to requiring that the reduced holonomy of
$\hat\nabla_W$ is included in $SU(n)$.

We shall not investigate the general case, instead we shall take
$M=X\times T^n$ with metric $g=g_X+g_T$ and Hermitian form
$\omega=\omega_X+\omega_T$ and assume that $X$ and $T^{2m}$
equipped with $(g_X, \omega_X)$ and $(g_T, \omega_T)$,
respectively, are Hermitian manifolds. Next we apply a Swann twist
associated with a $T^{2m}$ principal bundle over $M$ with
connection $\lambda$ and curvature $F$ supported on $X$. Then $W$
is a $T^{2m}$ fibration over $X$ with metric $g_W$ and Hermitian
form $\omega_W$ given by \bea\nonumber \pi_W^*g_W= h_{ab} \tilde
\lambda^a\otimes \tilde \lambda^b+ \pi^* g_X \cr
\pi_W^*\omega_W={1\over2}
J_{ab}\tilde\lambda^a\wedge\tilde\lambda^b+\pi^* \omega_X \eea
where $\tilde \lambda(\mathring \xi)=0$, $d\tilde \lambda=F$,
$h=g_T(\beta^{-1}\xi, \beta^{-1}\xi)$ and
$J=\omega_T(\beta^{-1}\xi, \beta^{-1}\xi)$.  In such a case,  one
can show, see\cite{grant, iibhor}, that $\hat\rho_W=0$, provided
that \bea\nonumber \hat\rho_X&=&-  \kappa^a h_{ab} F^b~,~~~ \cr
\omega_X\cdot F^a&=& \kappa^a~, \eea where $\hat\rho_X$ is the
Ricci form of the Hermitian connection with skew-symmetric torsion
on $X$, $\kappa$ is constant and $\omega_X\cdot F^a$ is the inner
product of $\omega_X$ and  $F$. Observe that if $F$ is
Hermitian-Einstein, ie $\kappa=0$, then $\hat\rho_X=0$ and so $X$
is CYT.

To find examples, let us take $M=X\times T^2$ with $X=X^{2k}\times X^2$, where both $X^{2k}\times X^2$ are K\"ahler-Einstein spaces with cosmological
constants $\ell_1$ and $\ell_2$, respectively. In such case, the Ricci forms of the K\"ahler metrics satisfy
 $\rho_{X^{2k}}=\ell_1 \omega_{(2k)}$ and $\rho_{X^{2}}=\ell_2 \omega_{(2)}$.
Assuming that $\omega_{(2k)} \in H^2(X^{2k}, \bZ)$ and
$\omega_{(2)} \in H^2(X^{2}, \bZ)$ and focusing on the k-SKT
examples for which $F^1=\omega_{(2k)}$ and $F^2=\omega_{(2)}$.
Then $\kappa=(k, 1)$, and so on finds that \bea\nonumber \ell_1=-k
h_{11}- h_{12}~,~~~\ell_2= -k h_{12}-h_{22}~. \eea
 Next consider the examples of k-SKT manifolds  with $g_T(\xi_A, \xi_B)=\delta_{AB}$ and $\beta$ given in (\ref{sktalpha}). Then one can show  that the above two conditions become
\bea\nonumber \ell_1&=&-k (p_1^2+p_2^2)- (p_1 q_1+p_2
q_2)=-{k+1\over2} (p_1^2+p_2^2)~,~~~ \cr \ell_2&=&-k (p_1 q_1+p_2
q_2)-(q_1^2+q_2^2)~, \eea where we have also used (\ref{ksktex2}).

Next consider  the k-SKT manifolds constructed from $M=X\times
T^2$ with $X=X^{2k}\times X^4$. Assuming  that both $X^{2k}$ and
$X^4$ are K\"ahler-Einstein, one finds that \bea\nonumber
\ell_1&=&-k (p_1^2+p_2^2)-2 (p_1 q_1+p_2 q_2)~, \cr \ell_2&=&-k
(p_1 q_1+p_2 q_2)-2 (q_1^2+q_2^2)~, \eea where $\ell_1$ and
$\ell_2$ are the cosmological constants of $X^{2k}$ and $X^4$,
respectively.

\vskip 0.5cm
\noindent{\bf Acknowledgements} \vskip 0.1cm
\noindent
SI is  partially supported by the Contract 181/2011
with the University of Sofia `St.Kl.Ohridski', Contract ``Idei", DID 02-39/21.12.2009 and  Contract ``Idei", DO 02-257/18.12.2008.
GP is partially supported by  the STFC rolling grant ST/G000/395/1.
\vskip 0.5cm


\end{document}